\documentclass[12pt]{article} 
\usepackage{latexsym} 
\usepackage{amsmath}
\usepackage{mathptmx}
\usepackage{amssymb} 
\usepackage{amsfonts}
\usepackage{amscd} 
\usepackage{graphicx}
\usepackage{layout}
\begin{document} 
\newtheorem{Th}{Theorem}[section]
\newtheorem{Cor}{Corollary}[section]
\newtheorem{Prop}{Proposition}[section]
\newtheorem{Lem}{Lemma}[section]
\newtheorem{Def}{Definition}[section]
\newtheorem{Rem}{Remark}[section]
\newtheorem{Ex}{Example}[section]
\newtheorem{stw}{Proposition}[section]


\newcommand{\bet}{\begin{Th}}
\newcommand{\ent}{\stepcounter{Cor}
   \stepcounter{Prop}\stepcounter{Lem}\stepcounter{Def}
   \stepcounter{Rem}\stepcounter{Ex}\end{Th}}


\newcommand{\bec}{\begin{Cor}}
\newcommand{\enc}{\stepcounter{Th}
   \stepcounter{Prop}\stepcounter{Lem}\stepcounter{Def}
   \stepcounter{Rem}\stepcounter{Ex}\end{Cor}}
\newcommand{\bep}{\begin{Prop}}
\newcommand{\enp}{\stepcounter{Th}
   \stepcounter{Cor}\stepcounter{Lem}\stepcounter{Def}
   \stepcounter{Rem}\stepcounter{Ex}\end{Prop}}
\newcommand{\bel}{\begin{Lem}}
\newcommand{\enl}{\stepcounter{Th}
   \stepcounter{Cor}\stepcounter{Prop}\stepcounter{Def}
   \stepcounter{Rem}\stepcounter{Ex}\end{Lem}}
\newcommand{\bef}{\begin{Def}}
\newcommand{\enf}{\stepcounter{Th}
   \stepcounter{Cor}\stepcounter{Prop}\stepcounter{Lem}
   \stepcounter{Rem}\stepcounter{Ex}\end{Def}}
\newcommand{\ber}{\begin{Rem}}
\newcommand{\enr}{
   \stepcounter{Th}\stepcounter{Cor}\stepcounter{Prop}
   \stepcounter{Lem}\stepcounter{Def}\stepcounter{Ex}\end{Rem}}
\newcommand{\bee}{\begin{Ex}}
\newcommand{\ene}{
   \stepcounter{Th}\stepcounter{Cor}\stepcounter{Prop}
   \stepcounter{Lem}\stepcounter{Def}\stepcounter{Rem}\end{Ex}}
\newcommand{\Proof}{\noindent{\it Proof\,}:\ }
\newcommand{\beP}{\Proof}
\newcommand{\enP}{\hfill $\Box$ \par\vspace{5truemm}}

\newcommand{\EE}{\mathbb{E}}
\newcommand{\QQ}{\mathbb{Q}}
\newcommand{\OO}{\mathbb{O}}
\newcommand{\R}{\mathbb{R}}
\newcommand{\C}{\mathbb{C}}
\newcommand{\ZZ}{\mathbb{Z}}
\newcommand{\KK}{\mathbb{K}}
\newcommand{\NN}{\mathbb{N}}
\newcommand{\PP}{\mathbb{P}}
\newcommand{\HH}{\mathbb{H}}
\newcommand{\uuu}{\boldsymbol{u}}
\newcommand{\xxx}{\boldsymbol{x}}
\newcommand{\aaa}{\boldsymbol{a}}
\newcommand{\bbb}{\boldsymbol{b}}
\newcommand{\AAA}{\mathbf{A}}
\newcommand{\BBB}{\mathbf{B}}
\newcommand{\ccc}{\boldsymbol{c}}
\newcommand{\iii}{\boldsymbol{i}}
\newcommand{\jjj}{\boldsymbol{j}}
\newcommand{\kkk}{\boldsymbol{k}}
\newcommand{\rrr}{\boldsymbol{r}}
\newcommand{\FFF}{\boldsymbol{F}}
\newcommand{\yyy}{\boldsymbol{y}}
\newcommand{\ppp}{\boldsymbol{p}}
\newcommand{\qqq}{\boldsymbol{q}}
\newcommand{\nnn}{\boldsymbol{n}}
\newcommand{\vvv}{\boldsymbol{v}}
\newcommand{\eee}{\boldsymbol{e}}
\newcommand{\fff}{\boldsymbol{f}}
\newcommand{\www}{\boldsymbol{w}}
\newcommand{\0}{\boldsymbol{0}}
\newcommand{\lon}{\longrightarrow}
\newcommand{\ga}{\gamma}
\newcommand{\pa}{\partial}
\newcommand{\QED}{\hfill $\Box$}
\newcommand{\id}{{\mbox {\rm id}}}
\newcommand{\Ker}{{\mbox {\rm Ker}}}
\newcommand{\grad}{{\mbox {\rm grad}}}
\newcommand{\ind}{{\mbox {\rm ind}}}
\newcommand{\rot}{{\mbox {\rm rot}}}
\newcommand{\diver}{{\mbox {\rm div}}}
\newcommand{\Gr}{{\mbox {\rm Gr}}}
\newcommand{\LG}{{\mbox {\rm LG}}}
\newcommand{\Diff}{{\mbox {\rm Diff}}}
\newcommand{\Symp}{{\mbox {\rm Symp}}}
\newcommand{\Ct}{{\mbox {\rm Ct}}}
\newcommand{\Uns}{{\mbox {\rm Uns}}}
\newcommand{\rank}{{\mbox {\rm rank}}}
\newcommand{\sign}{{\mbox {\rm sign}}}
\newcommand{\Spin}{{\mbox {\rm Spin}}}
\newcommand{\Sp}{{\mbox {\rm Sp}}}
\newcommand{\Int}{{\mbox {\rm Int}}}
\newcommand{\Hom}{{\mbox {\rm Hom}}}
\newcommand{\Tan}{{\mbox {\rm Tan}}}
\newcommand{\codim}{{\mbox {\rm codim}}}
\newcommand{\ord}{{\mbox {\rm ord}}}
\newcommand{\Iso}{{\mbox {\rm Iso}}}
\newcommand{\corank}{{\mbox {\rm corank}}}
\def\mod{{\mbox {\rm mod}}}
\newcommand{\pt}{{\mbox {\rm pt}}}
\newcommand{\qed}{\hfill $\Box$ \par}
\newcommand{\spe}{\vspace{0.4truecm}}
\renewcommand{\0}{\mathbf 0}
\newcommand{\ad}{{\mbox{\rm ad}}}
\newcommand{\xdownarrow}[1]{%
  {\left\downarrow\vbox to #1{}\right.\kern-\nulldelimiterspace}
}

\newcommand{\dint}[2]{{\displaystyle\int}_{{\hspace{-1.9truemm}}{#1}}^{#2}}

\title{Duality of (2,3,5)-distributions and Lagrangian cone structures
}

\author{Goo Ishikawa\thanks{Hokkaido University, Sapporo 060-0810, Japan. 
E-mail: ishikawa@math.sci.hokudai.ac.jp}, \, 
Yumiko Kitagawa\thanks{
Oita National College of Technology, Oita 870-0152, Japan. 
E-mail: kitagawa@oita-ct.ac.jp}, \, 
Asahi Tsuchida\thanks{Hokkaido University, Sapporo 060-0810, Japan. 
E-mail: asahi-t@math.sci.hokudai.ac.jp}, \, 
Wataru Yukuno\thanks{Hokkaido University, Sapporo 060-0810, Japan. 
E-mail: yukuwata@math.sci.hokudai.ac.jp}
}

\date{}

\maketitle

\begin{abstract}
As was shown by a part of the authors, for a given $(2, 3, 5)$-distribution $D$ on a $5$-dimensional manifold $Y$, there is, locally, a Lagrangian cone structure $C$ on another $5$-dimensional manifold $X$ which consists of abnormal or singular paths of $(Y, D)$.
We give a characterization of the class of Lagrangian cone structures corresponding to $(2, 3, 5)$-distributions. Thus we complete the duality between $(2, 3, 5)$-distributions and Lagrangian cone structures via pseudo-product structures of type $G_2$.
A local example of non-flat perturbations of the global model of flat Lagrangian cone structure which corresponds to $(2,3,5)$-distributions is given.
\end{abstract}

\section{Introduction}

A distribution $D$ on a $5$-dimensional manifold $Y$ is called a {\it $(2,3,5)$-distribution} 
if there is a local section $\eta_1, \eta_2$ of $D$ such that 
$$
\eta_1, \ \eta_2, \ [\eta_1, \eta_2], \ [\eta_1, \ [\eta_1, \eta_2]], \ [\eta_2, [\eta_1, \eta_2]]
$$
form a local frame of the tangent bundle to $Y$, in other words, if $D$ has the weak growth $(2, 3, 5)$, 
namely, if $\rank(\pa{\mathcal D}) = 3$ and $\rank(\pa^{(2)}{\mathcal D}) = 5$, 
where $\pa{\mathcal D} := 
[{\mathcal D}, {\mathcal D}] \ (= \mathcal D + [{\mathcal D}, {\mathcal D}])$, 
the derived system, and $\pa^{(2)}{\mathcal D} := 
[{\mathcal D}, \pa{\mathcal D}] \ (= \mathcal D + \pa{\mathcal D} + [{\mathcal D}, \pa{\mathcal D}])$ 
for the sheaf $\mathcal D$ of section-germs to $D$. 

The geometry and classification problem of $(2,3,5)$-distributions are studied after E. Cartan (\cite{Cartan10}),  related to the simple Lie group $G_2$, by many mathematicians
 (\cite{AZ}\cite{Bryant}\cite{LNS}\cite{The}\cite{Yamaguchi}\cite{Yamaguchi99}\cite{Zelenko}\cite{Zhitomirskii}). 
The $(2,3,5)$-distributions are related to many problems, for instance, 
to the problem of \lq\lq rolling balls" 
(\cite{Agrachev}\cite{BM}\cite{BH14}\cite{AN}), to indefinite conformal metrics 
(\cite{Nurowski05}\cite{LNS}), to non-linear differential equations 
(\cite{Randall}), and so on. 

In \cite{Kitagawa}\cite{IMT3}\cite{Kitagawa2}, we studied the global duality of $G_2$-homogeneous (flat) 
$(2, 3, 5)$-distribution and a Lagrangian cone structure 
from Cayley's split Octonions and classified the related generic singularities. 
In \cite{IKY}, we associated locally with any given $(2, 3, 5)$-distribution $D$ on a $5$-dimensional manifold $Y$, 
a Lagrangian cone structure $C$ on another $5$-dimensional manifold $X$, 
which consists of abnormal or singular paths of $(Y, D)$, in the sense 
of sub-Riemannian geometry or geometric control theory (see \cite{Montgomery}\cite{AS}). 
Moreover it was shown in \cite{IKY} that the original space $Y$ turns to be the totality of singular paths 
of the \lq\lq Lagrangian cone structure" $(X, C)$, 
when the cone field $C$ is regarded as a control system on $X$. 

In this paper, we give the characterization of the class of Lagrangian cone structures 
corresponding to $(2, 3, 5)$-distributions, 
and thus we complete the duality between $(2, 3, 5)$-distributions 
and Lagrangian cone structures (Theorem \ref{complete-duality}). 
The duality is actually understood via pseudo-product structure of $G_2$-type $E = K\oplus L$ 
on a $6$-dimensional manifold $Z$ (\S\ref{Pseudo-product structures}), 
which is regarded both as the prolongation of $(Y, D)$ 
and $(X, C)$ in the sense of Bryant (\cite{BH93}\cite{Bryant}), 
via the double fibration 
$$
(Y, D) \xleftarrow{\ \pi_Y\ } (Z, E) \xrightarrow{\ \pi_X\ } (X, C). 
$$

We realize, for the characterization a class of Lagrangian cone structures, that the language of cone structures 
is actually lacking: We introduce, regarding the cone structures as control systems, 
the notions of linear approximations and osculating bundles of cone structures, as well as the exact definition 
of non-degenerate Lagrangian cone structures (Definition \ref{Lagrangian-cone-structure-definition}). 

We remark that our correspondence is purely local in nature: It is \lq\lq spatially" local for $(Z, E)$ while 
\lq\lq spatially and directionally" local for $(Y, D)$ and for $(X, C)$. 
Moreover the \lq\lq directional locality"
for the distribution $(Y, D)$ is resolved by taking linear hull, 
however it is not the case for the cone structure $(X, C)$. This fact makes our duality delicate. 

It is clear that $(2,3,5)$-distributions form an open set, for Whitney $C^\infty$-topology, in the space of 
all distributions of rank $2$ on a $5$-dimensional manifold. In particular a $(2,3,5)$-distribution 
remains a $(2,3,5)$-distribution by sufficiently small perturbations with compact supports. 
However it is not clear such a stability for cone structures which corresponding to $(2, 3, 5)$-distributions 
via the duality. We give a local example of non-flat perturbations of the global model of 
flat Lagrangian cone structure (\cite{IMT3}), which corresponds to $(2,3,5)$-distributions (Example 
\ref{example-of-non-flat-Lagrange-cone-corresponding-to-(2,3,5)}). It is open 
the existence of non-flat {\it global} 
perturbations of Lagrangian cone structures which correspond to $(2, 3, 5)$-distributions. 

The cone structure was first given in \cite{AZ} by a foliation on the space
$P((\pa D)^\perp) \subset P(T^*Y)$ for the derived system $\pa D$ of a $(2, 3, 5)$-distribution $D$, 
which is an essentially same foliation in the space $P(D) \subset P(TY)$ of \cite{IKY}. See also \cite{AZ}\cite{DZ}. 
In fact there exists the natural fiber-preserving diffeomorphism $P(D) \to P((\pa D)^\perp)$ 
which preserves also the foliation induced from singular paths of $D$. 
Moreover the Lagrangian cone structure $C \subset TX$, which is contained in a contact structure 
$D' \subset TX$ on $X$,  has the essentially same information  
with the Jacobi curves introduced in \cite{AZ1}\cite{AZ2}. 
In fact each cone $C_x \subset D'_x, (x \in X)$ gives 
the (reduced) Jacobi curve associated to the singular path $x$ of $D$ 
in Lagrangian Grassmannian of $D'_x$ by taking tangent planes to $C_x$. 

In \cite{Zelenko}, it was shown that the Cartan tensor of any $(2, 3, 5)$-distribution is given by 
the fundamental invariant of Jacobi curves of singular paths and, in particular, the Cartan tensor is determined by 
the projective equivalence classes of the point-wise curves $P(C_x), x \in X$ 
of the corresponding Lagrangian cone structure $(X, C)$. 
We give a short proof (Proposition \ref{cubic-flat}), 
related to the study on $G_2$-contact structures (\cite{CS}\cite{LNS}), 
that the $(2, 3, 5)$-distribution which corresponds to a cubic Lagrangian cone structure via our duality 
is necessarily flat, by using Zelenko's theorem \cite{Zelenko}. 

In \S \ref{Pseudo-product structures}, 
we review the results given in the previous paper \cite{IKY} with additional explanations. In particular 
we give the exact definition of (non-degenerate) Lagrangian cone structures 
(Definition \ref{Lagrangian-cone-structure-definition}). 

In \S \ref{Complete duality}, we complete the duality between $(2,3,5)$-distributions and 
non-degenerate Lagrangian cone structures with an additional condition via 
pseudo-product structures of type $G_2$. 

We conclude this paper by several remarks related to the duality in \S \ref{cubic Lagrangian cone structures}. 

\

All manifolds and mappings are supposed to be of class $C^\infty$ unless otherwise stated. 

The authors are grateful to Professor Hajime Sato for valuable comment. 

\section{Pseudo-product structures of $G_2$-type}
\label{Pseudo-product structures}

Let $D$ be a $(2, 3, 5)$-distribution on a $5$-dimensional manifold $Y$. 
Let $Z := P(D) = (D - 0)/\R^{\times}$ be the space of tangential lines in $D$, 
$
Z := \{ (y, \ell) \mid y \in Y,  \ \ell \subset D_y (\subset T_yY), \, \dim(\ell) = 1\}. 
$
Then $\dim(Z) = 6$ and the projection $\pi_Y : Z \to Y$ is an $\R P^1$-bundle. 

We define a subbundle 
$E \subset TZ$ of rank $2$, {\it Cartan prolongation} of $D \subset TY$, by setting 
for each $(y, \ell) \in Z$, $\ell \subset D_y$, 
$E_{(y, \ell)} := \pi_{Y*}^{-1}(\ell) \ (\subset T_{(y, \ell)}Z)$. 
Then $E$ is a distribution with (weak) growth $(2, 3, 4, 5, 6)$: 
$\rank(E) = 2, \, \rank(\pa E) = 3, \, \rank(\pa^{(2)}E) = 4, \, $
$\rank(\pa^{(3)}E) = 5, \, \rank(\pa^{(4)}E) = 6. $

Then we see that there exists an {\it intrinsic} decomposition
$$
E = K \oplus L
$$ 
of $E$ with $L := \Ker(\pi_{Y*}) \subset E$ and 
a complementary line subbundle $K$ of $E$, 
a {\it pseudo-product structure} in the sense of N. Tanaka \cite{Tanaka79}\cite{Tanaka85}. 

We will explain this in terms of \lq\lq geometric control theory" (\cite{AS}\cite{Montgomery}). 

A {\it control system} ${\mathbb C} : {\mathcal U} \xrightarrow{F} TM \to M$ on a manifold 
$M$ is given by 
a locally trivial fibration $\pi_{\mathcal U} : {\mathcal U} \to M$ over $M$ and 
a map $F : {\mathcal U} \to TM$ such that the following diagram commutes: 
$$
\begin{array}{ccc}
{\mathcal{U}} \  & \ \xrightarrow{ \ F \ } &  \ TM
\vspace{0.2truecm}
\\
\ \ \ \ {\mbox{\footnotesize {$\pi_{\mathcal{U}}$}}}\searrow\hspace{-0.2truecm} &  
& \hspace{-0.2truecm}\swarrow{\mbox{\footnotesize {$\pi_{TM}$}}}  \ \ \ \
\vspace{0.2truecm}
\\
    & M &  
\end{array}
$$

Any section $s : M \to {\mathcal{U}}$ defines a vector field $F\circ c : M \to TM$ over $M$. 
Via a local triviality on $M$, a control system is given by 
a family of vector fields $f_u(x) = F(x, u)$ over $M$, 
$(x, u) \in {\mathcal U}, x \in M$. 

A distribution $D \subset TM$ is regarded as a control system 
${\mathbb D} : D \hookrightarrow TM \longrightarrow M$, 
by the inclusion. 

Two control systems 
${\mathbb C} : {\mathcal U} \xrightarrow{F} TM \xrightarrow{\pi_{TM}} M$ and 
${\mathbb C}' : {\mathcal U}' \xrightarrow{F'} TM' \xrightarrow{\pi_{TM'}} M'$ 
are called {\it isomorphic} if 
the diagram 
$$
\begin{array}{ccccc}
{\mathcal U} & \xrightarrow{\ F\ } & TM & \xrightarrow{\pi_{TM}} & M
\\
\psi\downarrow \ \ \ & & \varphi_*\!\downarrow\ \ \  & & \ \ \downarrow\varphi
\\
{\mathcal U}' & \xrightarrow{\ F'\ } & TM' & \xrightarrow{\pi_{TM'}} & M'
\end{array}
$$
commutes for some diffeomorphisms $\psi$ and $\varphi$. 
Here $\varphi_*$ is the differential of $\varphi$. 

The pair $(\psi, \varphi)$ of diffeomorphisms is called 
an {\it isomorphism} of the control systems ${\mathbb C}$ and 
${\mathbb C}'$. 

Given a control system ${\mathbb C} : {\mathcal U} \xrightarrow{F} TM \to M$, 
an $L^\infty$ (measurable, essentially bounded) 
map $c : [a, b] \to {\mathcal U}$ is called an {\it admissible control} 
if the curve 
$$
\gamma := \pi_{\mathcal U}\circ c : [a, b] \to M
$$ 
satisfies the differential equation 
$$
\dot{\gamma}(t) = F(c(t)) \quad ({\mbox{\rm a.e.}}\  t \in [a, b]). 
$$
Then the Lipschitz curve $\gamma$ is called a {\it trajectory}. 
If we write $c(t) = (x(t), u(t))$, then $x(t) = \gamma(t)$ and 
$$
\dot{x}(t) = F(x(t), u(t)), \quad ({\mbox{\rm a.e.}} \ t \in [a, b]). 
$$
We use the term \lq\lq {\it path}\rq\rq \ for a smooth ($C^\infty$)
immersive trajectory regarded up to parametrisation. 

The totality ${\mathcal C}$ 
of admissible controls 
$c : [a, b] \to {\mathcal U}$ with a given initial point $q_0 \in M$ 
is a Banach manifold. 
The {\it endpoint mapping} 
${\mbox{\rm End}} : {\mathcal C} \to M$ is defined by 
$$
{\mbox{\rm End}}(c) := \pi_{\mathcal U}\circ c(b).  
$$
An admissible control $c : [a, b] \to {\mathcal U}$ with the initial point 
$\pi_{\mathcal U}(c(a)) = q_0$  
is called {\it singular} or {\it abnormal}, 
if $c \in {\mathcal C}$ is a singular point of ${\mbox{\rm End}}$, namely if 
the differential ${\mbox{\rm End}}_* : T_c{\mathcal C} \to T_{{\mathcal E}(c)}M$ is not surjective.  
If $c$ is a singular control, 
then the trajectory $\gamma = \pi_{\mathcal U}\circ c$ is called a 
{\it singular trajectory} or an {\it abnormal extremal}. 

Let $D \subset TY$ be a $(2,3,5)$-distribution. 
Then, it can be shown that for any point $y$ of $Y$ 
and for any direction $\ell \subset D_y$,  
there exists uniquely a {\it singular} 
$D$-{\it path} (an immersed abnormal extremal for $D$) 
through $y$ with the given direction $\ell$. 
Thus the singular $D$-paths form another five dimensional manifold $X$. 

Let $Z = P(D) = (D - 0)/\R^{\times}$ be the space of tangential lines in 
$D$, $\dim(Z) = 6$. 
Then $Z$ is naturally foliated by the liftings of singular 
$D$-paths, and we have locally double fibrations: 
$$
Y \xleftarrow{\ \pi_Y\ } Z \xrightarrow{\ \pi_X\ } X. 
$$

If we put $L = \Ker(\pi_{Y*}), K = \Ker(\pi_{X*})$, 
then we have a decomposition $E = K \oplus L$ by sub-bundles of rank $1$. 

We denote, for any distribution $E$, by ${\mathcal E}$ the sheaf of local sections to $E$.  
We set 
$$
\pa{\mathcal E} := [{\mathcal E}, {\mathcal E}] = {\mathcal E} + [{\mathcal E}, {\mathcal E}]
, \quad \pa^{(2)}{\mathcal E} := [{\mathcal E}, \pa {\mathcal E}] = {\mathcal E} + \pa {\mathcal E} + [{\mathcal E}, \pa {\mathcal E}]
$$
and so on. 
If $\pa{\mathcal E}$ is generated by a local sections of a distribution, then we denote it by $\pa E$. 

\

\bef
{\rm 
A distribution $(Z, E)$ of rank $2$ on a $6$-dimensional manifold $Z$ with a decomposition 
$E = K \oplus L$ by subbundles $K, L$ of rank $1$ is called 
a pseudo-product structures of $G_2$-type if $E$ has small growth $(2,3,4,5,6)$ and moreover 
satisfies that 
$$
\left[ {\mathcal K}, {\mathcal L} \right] = \pa {\mathcal E}, \quad 
\left[ {\mathcal K}, \pa {\mathcal E} \right] = \pa^{(2)}{\mathcal E}, \ \left[ {\mathcal L}, \pa {\mathcal E} \right] = \pa {\mathcal E}, 
$$
$$
\left[\right. {\mathcal K}, \pa^{(2)} {\mathcal E} \left.\right] = \pa^{(3)}{\mathcal E}, \ \left[\right. {\mathcal L}, \pa^{(2)} {\mathcal E}\left.\right] = \pa^{(2)} {\mathcal E}, 
\quad 
\left[\right. {\mathcal K}, \pa^{(3)} {\mathcal E} \left.\right] = \pa^{(3)}{\mathcal E}, \ \left[\right. {\mathcal L}, \pa^{(3)} {\mathcal E}\left.\right] = \pa^{(4)} {\mathcal E}. 
$$
}
\enf

Then, by taking the gradation of the filtration
$$
{\mathcal E} \subset \pa{\mathcal E} \subset \pa^{(2)}{\mathcal E} \subset \pa^{(3)}{\mathcal E} \subset \pa^{(4)}{\mathcal E}, 
$$
we have, at each point $z \in Z$, the symbol algebra: 
$$
\begin{array}{c}
\mathfrak{m} = \mathfrak{g}_{-5}\oplus \mathfrak{g}_{-4}\oplus \mathfrak{g}_{-3}\oplus 
\mathfrak{g}_{-2}\oplus \mathfrak{g}_{-1} 
=\langle e_6 \rangle \oplus \langle e_5 \rangle \oplus \langle e_4 \rangle \oplus 
\langle e_3 \rangle \oplus \langle e_1, e_2 \rangle, 
\vspace{0.2truecm}
\\  
{[e_1, e_2] = e_3, \ [e_1, e_3] = e_4,\ [e_2, e_3] = 0, \ [e_1, e_4] = e_5, \ [e_2, e_4] = 0, \
[e_1, e_5] = 0, \ [e_2, e_5]=e_6, }
\end{array}
$$
with the decomposition $\mathfrak{g}_{-1} = {\mathfrak{k}}\oplus{\mathfrak{l}} = \langle e_1\rangle \oplus 
\langle e_2\rangle$. 

\

Then we have

\bet
There exists a natural bijective correspondence of local isomorphism classes 
between $(2, 3, 5)$-distributions and pseudo-product structures of $G_2$-type. 
\ent

\Proof
First let us make sure that 
the prolongation $E$ of a $(2, 3, 5)$-distribution $D$ on a $5$-dimensional manifold $Y$ 
has small growth $(2, 3, 4, 5, 6)$. 

Let $\eta_1, \eta_2$ be a local frame of $D$. 
Then, setting 
$$
\eta_3 := [\eta_1, \eta_2], \ \eta_4 := [\eta_1, \eta_3], \ \eta_5 := [\eta_2, \eta_3], 
$$
we have a a local frame $\eta_1, \eta_2, \eta_3, \eta_4, \eta_5$ of $TY$. 
For each $y \in Y$, directions in $D_y$ are, locally, parametrized via $\eta_1(y) + t\eta_2(y)$ 
$(t \in \R)$. 
Then, for any system of local coordinates 
$y = (y_1, y_2, y_3, y_4, y_5)$ of $Y$ centered at base point of $Y$, 
$(y, t)$ form a system of local coordinates of $Z$ 
such that $\pi_Y$ is expressed by $(y, t) \mapsto y$. 
We regard $\eta_1, \eta_2, \eta_3, \eta_4, \eta_5$ 
as vector-fields over $Z$.  Then 
$$
\zeta_1 := \eta_1 + t\eta_2 \quad \zeta_2 := \frac{\pa}{\pa t}, 
$$
form a local frame of $E$, and 
$\eta_1, \eta_2, \eta_3, \eta_4, \eta_5, \zeta_2$ 
of $TZ$. 

Since $[\zeta_1, \zeta_2] = [\eta_1 + t\eta_2, \zeta_2] = - \eta_2$, 
we have 
$$
\pa{\mathcal E} = \langle \zeta_1, \zeta_2, \eta_2\rangle = \langle \eta_1, \eta_2, \zeta_2\rangle, 
$$
which is of rank $3$. Here $\langle \zeta_1, \zeta_2, \eta_2\rangle$ means the distribution generated by 
$\zeta_1, \zeta_2, \eta_2$. 
Since $[\zeta_1, \eta_2] = [\eta_1 + t\eta_2, \eta_2] = \eta_3$ and $[\zeta_2, \eta_2] = 0$, we have 
$$
\pa^{(2)}{\mathcal E} = \langle\eta_1, \eta_2, \eta_3, \zeta_2\rangle, 
$$
which is of rank $4$. 
Since $[\zeta_1,\eta_3] = [\eta_1 + t\eta_2, \eta_3] = \eta_4 + t\eta_5$ and $[\zeta_2, \eta_3] = 0$, 
we have 
$$
\pa^{(3)}{\mathcal E} = \langle\eta_1, \eta_2, \eta_3, \eta_4 + t\eta_5, \zeta_2\rangle, 
$$
that is of rank $5$.  
Since $[\zeta_2, \eta_4 + t\eta_5] = \eta_5$,  we have $\pa^{(4)}{\mathcal E} = TZ$. 
Therefore $E$ has small growth $(2, 3, 4, 5, 6)$. 

Note that ${\mathcal L}$ is generated by $\zeta_2$. 
Moreover there exists a generator of ${\mathcal K}$ of form 
$\zeta_1 + e(y, t)\zeta_2$. 
In fact the function $e(y, t)$ is uniquely determined by 
the condition 
$\left[\right.{\mathcal K}, \pa^{(3)} {\mathcal E}\left.\right] = \pa^{(3)}{\mathcal E}$, 
which is equivalent to the condition 
$$
e\eta_5 + [\eta_1, \eta_4] + t[\eta_1, \eta_5] + t[\eta_2, \eta_4] + t^2[\eta_1, \eta_5] \equiv 0, \quad 
\mod.\  \pa^{(3)}{\mathcal E}. 
$$
Then other remaining conditions that $E = K \oplus L$ is a pseudo-product structure 
of type $G_2$ follow. 

Conversely suppose $E = K \oplus L$ is a pseudo-product structure of type 
$G_2$. 
Then $L$ is the Cauchy characteristic of $\pa E$ (see \cite{BCGGG}).  Let $Y$ be the leaf space of $L$, which is 
locally defined 
$5$ dimensional manifold. Moreover $Z$ has a system of local coordinates $(y, t)$ 
centered at the base point 
such that $\pi_Y$ is given by $(y, t) \mapsto y$. 
Let $D$ be the reduction of $\pa E$ by $L$. Take a local frame $\eta_1, \eta_2$ of $D$ such that, 
regarded as vector fields over $Z$, 
$\eta_1$ generates the quotient bundle $(\pa E)/E$. 
Moreover $\zeta_1 = \eta_1 + \varphi(y, t)\eta_2$ and $\zeta_2 = \pa/\pa t$ generates $K$ and $L$ 
respectively for some function $\varphi(y, t)$ with $\varphi(0,0) = 0$. 
Since 
$$
[\zeta_1, \zeta_2] = [\eta_1 + \varphi\eta_2, \zeta_2] = - (\pa\varphi/\pa t)\eta_2, 
$$
we have that $\pa\varphi/\pa t \not= 0$. Set $\zeta_3 := \eta_2$. Then 
$$
[\zeta_1, \zeta_3] = [\eta_1, \eta_2] + \eta_2(\varphi)\eta_2 \equiv [\eta_1, \eta_2] \quad  \mod. \pa {\mathcal E}. 
$$
Therefore $\eta_1, \eta_2, [\eta_1, \eta_2]$ are linearly independent point-wise on $Y$. 
We set $\zeta_4 := \eta_3 = [\eta_1, \eta_2]$ as a vector field over $Z$. 
Then 
$$
[\zeta_1, \zeta_4] 
= [\eta_1, \eta_3] + \varphi[\eta_2, \eta_3] -
\eta_3(\varphi)\eta_2 \equiv [\eta_1, \eta_3] + \varphi[\eta_2, \eta_3] \quad \mod. \pa^{(2)} {\mathcal E},  
$$
and $[\zeta_2, \zeta_4] = [\pa/\pa t, \eta_3] = 0$. 
Set $\eta_4 = [\eta_1, \eta_3], \eta_5 = [\eta_2, \eta_3]$ and $\zeta_5 = \eta_4 + \varphi\eta_5$. 
Then $\eta_4(0) \in (\pa^{(3)} E)_0 \setminus (\pa^{(2)} E)_0$. 
Then we have that $[\zeta_2, \zeta_5](0) \not\in (\pa^{(3)} E)_0$, while 
$
[\zeta_2, \zeta_5] = (\pa\varphi/\pa t)\eta_5(0). 
$
Therefore $\eta_5(0) \not\in (\pa^{(3)} E)_0$. Therefore $\eta_1, \eta_2, \eta_3, \eta_4, \eta_5$ 
are linearly independent point-wise. Thus we see that $D$ is a $(2, 3, 5)$-distribution. 

These correspondences induce the bijection between local isomorphism classes of 
$(2,3,5)$-distributions and pseudo-product structures of $G_2$-type on $5$-manifolds. 
\qed

\

Note that the original $(2, 3, 5)$-distribution $D$ is obtained as the linear hull of 
the cone field (\lq\lq bowtie") induced from $K$: 
$$
D_y = {\mbox{\rm linear hull}}\left(\bigcup_{z \in \pi_Y^{-1}(y)}\pi_{Y*}(K_z) \subset T_yY\right).
$$
Also, the $(2, 3, 5)$-distribution $D$ is obtained as the reduction of 
$\pa E$ by Cauchy characteristic $L = \Ker(\pi_{Y*})$. 
 
On the other hand we obtain a cone field $C \subset TX$ on $X$ by setting, 
for each $x \in X$, 
$$
C_x := \bigcup_{z \in \pi_X^{-1}(x)}\pi_{X*}(L_z) \subset T_xX. 
$$

Now, to make sure, we formulate exactly the notion of Lagrangian cone structures (see \cite{CS}): 

\bef
\label{Lagrangian-cone-structure-definition}
{\rm 
{\rm (1)} Let $X$ be a manifold of dimension $m$. A subset $C \subset TX$ is called a {\it cone structure} if 
there is an $\R^{\times}$-invariant subset $\underline{C} \subset \R^m$, a model cone,  
such that, for any $x \in X$, there exists an open neighborhood $U$ of $x$ and a local triviality 
$\Phi : \pi^{-1}(U) \to U\times \R^m$ of $\pi : TX \to X$ over $U$ satisfying 
$\Phi(\pi^{-1}(U) \cap C) = U \times \underline{C}$. 

{\rm (2)} Suppose that the model cone 
$\underline{C}$ is non-singular away from the origin in $\R^m$. Then $P(C)$ is a submanifold of $P(TX)$. 
For each section $s : X \to P(C)$ for the projection $P(C) \to X$, 
we have the subbundle $T_sC \subset TX$ by taking tangent planes 
of $C_x$ along the direction $s(x)$ at every point $x \in X$. 
We call the distribution $T_sC$ the {\it linear approximation} of $C$ along $s$. 

{\rm (3)} A cone structure $C \subset TX$ is called a {\it Lagrangian cone structure} if  
there exists a contact structure $D'$ on $X$ such that 
$C \subset D'$ and, for any section $s : X \to P(C)$, $T_sC$ is a Lagrangian subbundle 
of $D'$. The last condition is equivalent to that, for any $x \in X$, 
$C_x \setminus \{ 0\}$ is a Lagrangian submanifold of the linear symplectic manifold 
$D_x'$, or equivalently, 
$P(C_x)$ is a Legendrian submanifold of the contact manifold $P(D'_x)$ induced from 
the conformal symplectic vector space $D_x'$. 

{\rm (4)} Let $\dim(X) = 5$. A Lagrangian cone structure $C \subset TX$ for a contact structure 
$D' \subset TX$ is 
called {\it non-degenerate} if the spatial projective curve segment $P(C_x) \subset P(D_x') \cong P^3$ 
is non-degenerate, i.e. the first, second and third derivatives of a parametrization 
of $P(C_x)$ are linearly independent. 
}
\enf

From the condition (4), for each direction field $s$ of $C$, 
we define {\it osculating bundles} 
$O^{(2)}_sC \subset TX$ of rank $3$ and $O^{(3)}_sC \subset TX$ of rank $4$, 
generated by osculating planes $O_2$ and $3$-dimensional osculating spaces $O_3$ to $P(C_x)$ with direction $s$. 
Then the contact structure $D'$ coincides with $O^{(3)}_sC$ which is independent of $s$. 

\

Because distributions are regarded as cone structures of special type,
the notion of Lagrangian cone structures is a natural generalization for that of Lagrangian 
subbundle of the tangent bundle over a contact manifold. 

\bel
In our case, the above $C \subset TX$ corresponding to a $(2, 3, 5)$-distribution $D \subset TY$ 
is a non-degenerate Lagrangian cone structure in the sense of Definition 
\ref{Lagrangian-cone-structure-definition}. 
\enl

\Proof
By the condition
$\left[\right.{\mathcal K}, \pa {\mathcal E}\left.\right] = \pa^{(2)}{\mathcal E}$, 
$C$ satisfies the conditions (1)(2) of Definition \ref{Lagrangian-cone-structure-definition}. 
By the condition $\left[\right.{\mathcal K}, \pa^{(3)} {\mathcal E}\left.\right] = \pa^{(3)}{\mathcal E}$, 
$K$ is the Cauchy characteristic of $\pa^{(3)}{\mathcal E}$. Then 
the distribution $D' \subset TX$ induced from $\pa^{(3)}{\mathcal E}$ is a contact structure by
the condition $\left[\right.{\mathcal L}, \pa^{(3)} {\mathcal E}\left.\right] = \pa^{(4)} {\mathcal E}$. 
Moreover $\pa^{(2)}{\mathcal E}$ projects to tangent spaces to $C_x$ 
along $\pi_X^{-1}(x)$. For any section $s : X \to L$, $s(x) \not= 0$, we have 
that the linear approximation $T_sC$ is a Lagrangian subbundle of $D'$ by the condition 
$\left[\right.{\mathcal L}, \pa^{(3)} {\mathcal E}\left.\right] = \pa^{(3)} {\mathcal E}$. 
Therefore $C$ satisfies also the condition (3) of Definition \ref{Lagrangian-cone-structure-definition}. 
Thus $(X, C)$ is a Lagrangian cone structure. 
Moreover by the condition 
$\left[\right.{\mathcal K}, \pa^{(2)} {\mathcal E}\left.\right] = \pa^{(3)}{\mathcal E}$, 
the condition (4) of Definition \ref{Lagrangian-cone-structure-definition} is satisfied. 
Therefore $(X, C)$ is a non-degenerate Lagrangian cone structure. 
\QED

\

Now, we regard the cone field $C \subset TX$ as a control system over $X$: 
$$
{\mathbb C} : L  \xrightarrow{\pi_{X*}\vert_L} TX \to X.
$$ 

Then we have showed in \cite{IKY} the following theorem: 

\bet
{\rm (Duality Theorem \cite{IKY})} 
Singular paths of the control system 
$$
{\mathbb C} : L  \xrightarrow{\pi_{X*}\vert_L} TX \to X
$$ 
are given by $\pi_X$-images of $\pi_Y$-fibers. 

Therefore, for any $x \in X$ and for any direction 
$\ell \subset C_x$, there exists uniquely a singular ${\mathbb C}$-paths 
passing through $x$ with the direction $\ell$ at $x$. 

Thus the original space $Y$ is identified with 
the space of singular paths for $(X, C)$, while $X$ is 
the space of singular paths for $(Y, D)$. 
\ent

\

We recall the local characterization of singular controls. 

For a control system 
${\mathbb C} : {\mathcal U} \xrightarrow{F} TM \to M$ on a manifold 
$M$, we consider the fibre-product ${\mathcal U} \times_MT^*M$, and 
define the {\it Hamiltonian function} 
$H : {\mathcal U} \times_MT^*M \to \R$ of the control system 
$F : {\mathcal U} \to TM$ by 
$$
H(x, p, u) := \langle p, F(x, u)\rangle, \quad ((x,u), (x,p)) \in {\mathcal U} \times_MT^*M. 
$$

A singular control $(x(t), u(t))$ is characterized by the liftability to 
an {\it abnormal bi-extremal} $(x(t), p(t), u(t))$ satisfying the constrained Hamiltonian equation 
$$
\begin{cases}
\ \dot{x}_i(t)  =  \ \ \dfrac{\pa H}{\pa p_i}(x(t), p(t), u(t)), \quad (1 \leq i \leq m)
\vspace{0.2truecm}
\\
\ \dot{p}_i(t)  =  - \dfrac{\pa H}{\pa x_i}(x(t), p(t), u(t)), \quad (1 \leq i \leq m)
\vspace{0.2truecm}
\\
\ \dfrac{\pa H}{\pa u_j}(x(t), p(t), u(t)) = 0, \quad (1 \leq j \leq r), \qquad 
p(t) \not= 0. 
\end{cases}
$$

Let $E \subset TZ$ be a distribution on a manifold $Z$ regarded as a control system. 
A singular path $x(t)$ for $E \subset TZ$ is called {\it regular singular} 
if it is associated with an abnormal bi-extremal 
$(x(t), p(t), u(t))$ such that $p(t) \in (\pa E)^\perp \setminus (\pa^{(2)}E)^\perp \subset T^*Z$. 
A singular path $x(t)$ for $E \subset TZ$ is called {\it totally irregular singular} 
if any associated abnormal bi-extremals $(x(t), p(t), u(t))$ satisfies that 
$p(t) \in (\pa^{(2)}E)^\perp \subset T^*Z$. 

From the pseudo-product structure on $E \subset TZ$, we have 

\bet
{\rm (Asymmetry Theorem \cite{IKY})} 
A singular path for $E \hookrightarrow TZ \to Z$ is either a $\pi_Y$-fibre 
or a $\pi_X$-fibre. Each $\pi_Y$-fibre is regular singular, while 
each $\pi_X$-fibre is totally irregular singular. 
\ent

\section{Complete duality}
\label{Complete duality}

The description of the duality on $(2, 3, 5)$-distributions $(Y, D)$ 
and non-degenerate Lagrangian cone structures $(X, C)$ via $(Z, E)$ 
which is given in \S \ref{Pseudo-product structures} 
should be completed 
by answering the question: \ 
{\it 
What kinds of non-degenerate 
Lagrangian cone structures do they correspond to $(2, 3, 5)$-distributions ?  
}

Then we have 

\bet
\label{complete-duality}
There exist natural bijective correspondences of isomorphism classes: 
\\
$\{ (2,3,5){\mbox{\rm -distributions}}\  (Y, D) \}\!/\!\cong$ 
$\longleftrightarrow$
$
\left\{
\begin{array}{c}
{\mbox{\rm pseudo-product structures of }} G_2 {\mbox{\rm -type}} \ $(Z, E)$: 
\\
(2,3,4,5,6){\mbox{\rm -distributions}}\  E \ {\mbox{\rm with a decomposition}}
\\
E =  K \oplus L, \ 
\rank(K) = \rank(L) = 1, 
\\
\left[ {\mathcal K}, {\mathcal L} \right] = \pa {\mathcal E} \ (:= [{\mathcal E}, {\mathcal E}] = {\mathcal E} + [{\mathcal E}, {\mathcal E}]), 
\\
\left[\right.{\mathcal K}, \pa {\mathcal E}\left.\right] = \pa^{(2)}{\mathcal E}, \ \left[\right.{\mathcal L}, \pa {\mathcal E}\left.\right] = \pa {\mathcal E}, 
\\
\left[\right.{\mathcal K}, \pa^{(2)} {\mathcal E}\left.\right] = \pa^{(3)}{\mathcal E}, \ \left[\right.{\mathcal L}, \pa^{(2)} {\mathcal E}\left.\right] = \pa^{(2)} {\mathcal E}, 
\\
\left[\right.{\mathcal K}, \pa^{(3)} {\mathcal E}\left.\right] = \pa^{(3)}{\mathcal E}, \ \left[\right.{\mathcal L}, \pa^{(3)} {\mathcal E}\left.\right] = \pa^{(4)} {\mathcal E}. 
\\
\end{array}
\right\}\!/\!\cong
$
\\
\hspace{5.6truecm}
$\longleftrightarrow$ 
$\left\{ 
\begin{array}{cc}
{\mbox{\rm non-degenerate Lagrangian cone structures}} \ (X, C) 
\\
{\mbox{\rm on $5$-dimensional manifolds $X$ with the condition}}
\\
\pa(T_sC) \subset O^{(2)}_sC, \ {\mbox{\rm for any direction field\ }} s \ {\mbox{\rm of\ }} C.
\end{array}
\right\}\!/\!\cong$
\ent

\

\noindent
{\it Proof of Theorem \ref{complete-duality}.} 

Let $X$ be a $5$-dimensional manifold 
and 
$C \subset TX$ a non-degenerate Lagrangian cone structure 
(Definition \ref{Lagrangian-cone-structure-definition}). 
Then $Z = P(C) := (C \setminus ({\mbox{\rm zero-section}}))/\R^\times$ is 
a $6$-dimensional manifold and that $\pi_X : Z \to X$ is 
a $C^\infty$-fibration with projective curves 
$P(C_x) \subset P(T_xX) \cong P^4$ as fibers. 

By the non-degeneracy condition, we have
that the first, second and third derivatives are linearly independent everywhere on 
$P(C_x)$, for any $x \in X$. 

Then we define a subbundle $E \subset TZ$ of rank $2$ by setting 
$$
E_{(x, \ell)} := (\pi_X)_*^{-1}(\ell), 
$$
for each $(x, \ell) \in Z$ as the prolongation of the cone structure $C \in TX$. 
We set $K = \Ker((\pi_X)_*)$. 

Let $x = (x_1, x_2, x_3, x_4, x_5)$ be a system of local coordinates of $X$ 
and $x, \theta$ that of $Z$ such that 
$\pi_X : Z \to X$ is given by $(x, \theta) \mapsto x$ and 
$E$ is generated by $\zeta_1 = \frac{\pa}{\pa\theta}$ and a vector field $\zeta_2(x, \theta)$ 
of form 
$$
\zeta_2(x, \theta) = 
\frac{\pa}{\pa x_1} + A\frac{\pa}{\pa x_2} + B\frac{\pa}{\pa x_3} + S\frac{\pa}{\pa x_4} + T\frac{\pa}{\pa x_5}, 
$$
where $A, B, S, T$ are function-germs of $x, \theta$. 
The projective curve $C_x \subset P(T_xX)$ is given by 
$$
\theta \mapsto [ 1 : A(x, \theta) : B(x, \theta) : S(x, \theta) : T(x, \theta) ]
$$
in homogeneous coordinates, for each $x \in X$. 

We have, on $Z$, 
$$
[\zeta_1, \zeta_2](x, \theta) =  \frac{\pa\zeta_2}{\pa\theta}(x, \theta) =: \zeta_3, 
$$
and 
$$
[\zeta_1, \zeta_3](x, \theta) = \frac{\pa^2\zeta_2}{\pa\theta^2}(x, \theta) =: \zeta_4. 
$$
In local coordinates, 
$$
\zeta_3 = A_\theta\frac{\pa}{\pa x_2} + B_\theta\frac{\pa}{\pa x_3} + S_\theta\frac{\pa}{\pa x_4} + T_\theta\frac{\pa}{\pa x_5}, 
\quad 
\zeta_4 = A_{\theta\theta}\frac{\pa}{\pa x_2} + B_{\theta\theta}\frac{\pa}{\pa x_3} + S_{\theta\theta}\frac{\pa}{\pa x_4} + T_{\theta\theta}\frac{\pa}{\pa x_5}, 
$$
and 
$$
\zeta_5 = A_{\theta\theta\theta}\frac{\pa}{\pa x_2} + B_{\theta\theta\theta}\frac{\pa}{\pa x_3} + S_{\theta\theta\theta}\frac{\pa}{\pa x_4} + T_{\theta\theta\theta}\frac{\pa}{\pa x_5}. 
$$

Any direction field $s$ of $C$ is given by $x \mapsto (x, \theta(x))$ for some functions $\theta(x)$ on $x$
and the linear approximation $T_sC$ of $C$ along a direction field $s$ is generated by 
by $\zeta_2(x, \theta(x)), \frac{\pa\zeta_2}{\pa\theta}(x, \theta(x))$. 
Moreover the osculating bundles $O^{(2)}_sC$ and $O^{(3)}_sC$ are generated by 
$\zeta_2(x, \theta(x)), \frac{\pa\zeta_2}{\pa\theta}(x, \theta(x)), \frac{\pa^2\zeta_2}{\pa\theta^2}(x, \theta(x))$ 
and by $\zeta_2(x, \theta(x)), \frac{\pa\zeta_2}{\pa\theta}(x, \theta(x)), \frac{\pa^2\zeta_2}{\pa\theta^2}(x, \theta(x)), 
\frac{\pa^3\zeta_2}{\pa\theta^3}(x, \theta(x))$ respectively. 

By the condition $\pa(T_sC) \subset O^{(2)}_sC$, we have that 
$[\zeta_2, \zeta_3] \equiv 0, \mod. \  \zeta_1, \zeta_2, \zeta_3, \zeta_4$. 
Then there exists uniquely a function $U(x, \theta)$ such that 
$\widetilde{\zeta}_2 = \zeta_2 + U\zeta_1$ is a Cauchy characteristic vector field of $\pa E$, so that 
$[\widetilde{\zeta}_2, \zeta_3] \equiv 0, \ \mod. \ \zeta_1, \zeta_2, \zeta_3$. 

Taking the subbundle $L \subset E$ generated by $\widetilde{\zeta}_2$, 
we have a pseudo-product structure $E = K \oplus L$ on $Z$ satisfying the conditions 
$$
\left[ {\mathcal K}, {\mathcal L} \right] = \pa {\mathcal E}, \quad 
\left[ {\mathcal K}, \pa {\mathcal E} \right] = \pa^{(2)}{\mathcal E}, 
\quad \left[ {\mathcal L}, \pa {\mathcal E} \right] = \pa {\mathcal E}, 
\quad \left[\right. {\mathcal K}, \pa^{(2)} {\mathcal E} \left.\right] = \pa^{(3)}{\mathcal E}. 
$$
Since we have, by Jacobi identity, 
$
[\widetilde{\zeta}_2, [\zeta_1, \zeta_3]] + [\zeta_1, [\zeta_3, \widetilde{\zeta}_2] + [\zeta_3, [\widetilde{\zeta}_2, \zeta_1]] = 0, 
$
we have that 
$$
[\widetilde{\zeta}_2, \zeta_4] 
\equiv [\zeta_1, [\widetilde{\zeta}_2, \zeta_3]] \equiv 0, \ \mod. \ \zeta_1, \zeta_2, \zeta_3, \zeta_4. 
$$
Therefore the condition $\left[\right. {\mathcal L}, \pa^{(2)} {\mathcal E}\left.\right] = \pa^{(2)} {\mathcal E}$ 
is satisfied. Because $O^{(3)}_sC \subset TX$ 
is independent of $s$ and is a contact structure on $X$, we have that 
$\left[\right. {\mathcal K}, \pa^{(3)} {\mathcal E} \left.\right] = \pa^{(3)}{\mathcal E}$ 
and that $\left[\right. {\mathcal L}, \pa^{(3)} {\mathcal E}\left.\right]$ generates the total tangent bundle 
$TZ$, therefore, the last condition 
$\left[\right. {\mathcal L}, \pa^{(3)} {\mathcal E}\left.\right] = \pa^{(4)} {\mathcal E}$. 

Thus we see that, if $C$ is a non-degenerate Lagrangian cone structure with the condition that 
$\pa(T_sC) \subset O^{(2)}_sC$ for any direction field $s$ of $C$, 
then $E = K \oplus L$ is a pseudo-product structure of $G_2$-type. 

This completes the proof of Theorem \ref{complete-duality}. 
\QED

\ber
{\rm 
The cone structure $C \subset TX$ is regarded as the control system over $X$, 
$$
{\mathbb C} : L \to TX \to X, 
\quad 
L \ni ((x, \ell), v) \mapsto (x, v) \mapsto x, 
$$
with $2$-control parameters. 
In local coordinates, the control system ${\mathbb C}$ is given by 
$$
F(x; r, \theta) := r\left( 
\frac{\pa}{\pa x_1} + A(x, \theta)\frac{\pa}{\pa x_2} + B(x, \theta)\frac{\pa}{\pa x_3} + S(x, \theta)\frac{\pa}{\pa x_4} + T(x, \theta)\frac{\pa}{\pa x_5}\right), 
$$
with the control parameters $r, \theta$. 
}
\enr

\section{$(2,3,5)$-distributions and cubic Lagrangian cone structures}
\label{cubic Lagrangian cone structures}

Let us denote by $G_2'$ the automorphism group of the split octonion algebra ${\mathbb O}'$. 
Then for a Borel group subgroup $B$ and parabolic subgroups $P_1, P_2$ containing $B$ of $G_2'$, 
we have a double fibration 
$$
Y = G_2'/P_1 \xleftarrow{\ \pi_Y\ } Z = G_2'/B \xrightarrow{\ \pi_X\ } X = G_2'/P_2, 
$$
a $(2,3,5)$-distribution $D \subset TY$ on $Y$, a pseudo-product structure of type $G_2$ $E = K\oplus L 
\subset TZ$ 
on $Z$ and a non-degenerate Lagrangian cubic cone structure $C \subset TX$ (see {IMT3}). 
It is known also that $Y$ is diffeomorphic to $S^3\times S^2$ (resp. $Z$ to 
$S^3\times S^3$, $X$ to $S^2\times S^3$). 
On each of three places, there exists Cartan's parabolic geometry as a natural 
non-flat geometry modeled on the homogeneous space. 
On $Y$ it is the geometry of $(2,3,5)$-distributions. On $Z$ it is the geometry of 
pseudo-product structures of type $G_2$. 
On $X$ it is $G_2$-contact structures (\cite{CS}\cite{LNS}). 
Moreover any $G_2$-contact structure is accompanied with and is recovered from 
a non-degenerate Lagrangian cubic cone structure. 

Hajime Sato \cite{Sato} has suggested to the first author that 
any $G_2$-contact structure corresponding to a $(2,3,5)$-distribution should be flat, 
from the exact comparison of curvatures for associated Cartan connections 
on pseudo-product $G_2$-structure 
and on $G_2$-contact structures (\cite{Tanaka79}\cite{Yamaguchi}). 
Here we would like to provide alternative proof for the fact. In fact we have: 

\bep
\label{cubic-flat}
Any $(2,3,5)$-distribution $(Y, D)$ which corresponds to a cubic cone structure $(X, C)$ must 
be flat. 
\enp

\noindent
{\it Proof of Proposition \ref{cubic-flat}.} 
For each $x \in X$, the cone 
$C_x \subset D'_x (\subset T_xX)$ gives 
the (reduced) \lq\lq Jacobi curve" in the sense of 
Agrachev and Zelenko \cite{AZ1}\cite{AZ2}\cite{Zelenko}\cite{AZ}. 
Then, in \cite{Zelenko}, it is proved that \lq\lq Cartan tensor" of $D$ is recovered by a projective 
invariant, the fundamental invariant, a kind of cross ratio, of $P(C_x)$ point-wise. 
In fact, for the cone $C_x \subset D_x \cong \R^4$, 
there is associated a curve $P(C_x)$ in Grassmannian $\Gr(2, \R^4)$, and 
the fundamental invariants is calculated from $P(C_x)$ in projective invariant way. 

Suppose a $(2,3,5)$-distribution $D$ 
corresponds to a cubic cone structure $C \subset D' \subset TX$. 
Then the cone structure is non-degenerate. 
Since all non-degenerate cubic cones are projectively equivalent point-wise, 
the Cartan tensor of $D$ coincides with the flat $(2,3,5)$-distribution. 
Therefore $D$ must be flat. 
\QED

\bee
{\rm (Cubic Lagrangian cone structures not corresponding to $(2,3,5)$-distributions.)
\\
Consider the cubic cone structure $C$ on $(\R^5, 0)$ around the direction $\theta = 0$, 
$$
\begin{array}{l}
F(x; r, \theta) = r\left(
\frac{\pa}{\pa x_1} + \theta\frac{\pa}{\pa x_2} + (\theta^2 + a)\frac{\pa}{\pa x_3} \right.
+ (\theta^3 - 3\theta a)\frac{\pa}{\pa x_4} 
\vspace{0.1truecm}
\\
\left. 
\hspace{8truecm} 
+ \{ x_3\theta -2x_2(\theta^2 + a) + x_1(\theta^3 - 3\theta a)\}\frac{\pa}{\pa x_5}
\right), 
\end{array}
$$
defined by a $C^\infty$ function $a(x_1)$ 
with $a(0) = 0$. 

Then $C$ is a non-degenerate Lagrangian cone structure for 
the contact structure $D' : dx_5 - x_3 dx_2 + 2x_2 dx_3 - x_1 dx_4 = 0$. 
Moreover 
$C$ satisfies the condition $\pa(T_sC) \subset O^{(2)}_sC$ for any $s : X \to L \setminus\{ 0\}$, 
 to correspond to a $(2, 3, 5)$-distribution, 
if and only if 
$a \not\equiv 0$. 
The case $a \equiv 0$ corresponds to the $G_2$-homogeneous flat case (\cite{IMT3}). 
}
\ene

\

The following gives examples of non-degenerate Lagrangian non-cubic cone structures which 
correspond to $(2, 3, 5)$-distributions and shows the necessity of the additional condition 
$\pa(T_sC) \subset O^{(2)}_sC$ of Theorem \ref{complete-duality}. 

\bee
\label{example-of-non-flat-Lagrange-cone-corresponding-to-(2,3,5)}
{\rm (Non-cubic Lagrangian cone structures corresponding to $(2,3,5)$-distributions.) 
\\
Consider a cone structure on $(\R^5, 0)$ around the direction $\theta = 0$, 
$$
\begin{array}{l}
F(x; r, \theta) = r\left(
\frac{\pa}{\pa x_1} + \theta\frac{\pa}{\pa x_2} + (\theta^2 + b)\frac{\pa}{\pa x_3} 
+ (\theta^3 + c)\frac{\pa}{\pa x_4} \right.
\\
\hspace{8truecm}
\left. + \{ x_3\theta -2x_2(\theta^2 + b) + x_1(\theta^3 + c)\} \frac{\pa}{\pa x_5}
\right), 
\end{array}
\vspace{-0.3truecm}
$$
where $b = b(\theta), c = c(\theta)$, with $\ord_0 b(\theta) \geq 3, \ \ord_0 c(\theta) \geq 4$. 

Then $F$ is a non-degenerate Lagrangian cone structure,  
for the contact structure $D' : dx_5 - x_3 dx_2 + 2x_2 dx_3 - x_1 dx_4 = 0$. 
Moreover $F$ satisfies 
the condition $\pa(T_sC) \subset O^{(2)}_sC$, for any direction field $s$,  
to correspond to a $(2, 3, 5)$-distribution, 
if and only if 
$
c_\theta = 3\theta b_\theta - 3b. 
$

If $b_{\theta\theta\theta\theta} \not= 0$, then $C$ is not cubic. 
}
\ene

\end{document}